\documentclass[12pt]{amsart}
\usepackage{a4}
\usepackage{amssymb}

\newtheorem{thm}{Theorem}[section]

\newtheorem{prop}[thm]{Proposition}

\theoremstyle{definition}

\newtheorem{example}[thm]{Example}
\theoremstyle{remark}

\numberwithin{equation}{section}

\begin{document}

\title{Randomized series and Geometry of Banach spaces}

\author{Han Ju Lee}
\address{Mathematics Department 202 Mathematical Sciences Bldg University of
Missouri Columbia, MO 65211 USA } \email{hahnju@postech.ac.kr}

\thanks{This work was supported by grant  No. R01-2004-000-10055-0 from the
Basic Research Program of the Korea Science \& Engineering
Foundation}








\begin{abstract} We study some properties of the randomized series and their
applications to the geometric structure of Banach spaces. For $n\ge
2$ and $1<p<\infty$, it is shown that $\ell_\infty^n$ is
representable in a Banach space $X$ if and only if it is
representable in the Lebesgue-Bochner $L_p(X)$. New criteria for
various convexity properties in Banach spaces are also studied. It
is proved that a Banach lattice $E$ is uniformly monotone if and
only if its $p$-convexification $E^{(p)}$ is uniformly convex and
that a K\"othe function space $E$ is upper locally uniformly
monotone if and only if its $p$-convexification $E^{(p)}$ is
midpoint locally uniformly convex.
\end{abstract}

\maketitle

\section{Randomized series}
Let $\{r_n\}_{n=1}^\infty$ be a sequence of mutually independent,
symmetric  and integrable random variables on a probability space
$(\Omega, \mathcal{F}, P)$ and $\{x_n\}_{n=0}^\infty$ an arbitrary
sequence in a Banach space $X$. A {\it randomized series} $S_n$  is
a vector-valued random variable defined by
\[ S_n = x_0 + r_1 x_1 + \cdots + r_n x_n,\ \ \ (n=0, 1,\ldots)\]

Let $\mathcal{F}_0$ be the trivial $\sigma$-algebra $\{\emptyset,
\Omega\}$ and $\mathcal{F}_k$, $(k\ge 1)$  the $\sigma$-algebra
generated by random variables $\{r_i\}_{i=1}^k$. It is easy to see
that the sequence of randomized series $\{S_n\}_{n=0}^\infty$ is a
martingale with respect to the filtration
$\{\mathcal{F}_n\}_{n=0}^\infty$. In this paper, $\mathbb{E}$ stands
for the expectation with respect to the probability $P$.

We begin with the basic properties of the randomized series. For
more properties of random series, see \cite{K}.
\begin{prop}\label{propsubmartingale}
Let $\varphi$ be a convex function on $\mathbb{R}$. Then
$\{\varphi(\|S_n\|)\}_{n=0}^\infty$ is a submartingale with respect
to $\{\mathcal{F}_{n}\}^\infty_{n=0}$. In particular,
\[\mathbb{E}[\varphi(\|S_0\|)]\le \mathbb{E}[\varphi(\|S_1\|)]\le
\mathbb{E}[\varphi(\|S_2\|)]\le \cdots.\]
\end{prop}
\begin{proof}In the proof, the notation $\mathbb{E}[\;\cdot\;
|\mathcal{G}]$ means the conditional expectation with respect to the
sub-$\sigma$-algebra $\mathcal{G}$ and we need to show that
$\mathbb{E}[\varphi(\|S_{n+1}\|)|\mathcal{F}_n]\ge \varphi(\|S_n\|)$
almost surely. By the independence and symmetry of $\{r_k\}$ we get,
for almost all $\omega\in \Omega$, \begin{align*}
\mathbb{E}[\varphi&(\|S_{n+1}\|)|\mathcal{F}_n](\omega) =
\int_{\Omega} \varphi(\|S_n(\omega) + r_{n+1}(t)x_{n+1}\|) \;dP(t)\\
&= \int_{\Omega} \frac{\varphi(\|S_n(\omega) + r_{n+1}(t)x_{n+1}\|)
+ \varphi(\|S_n(\omega) - r_{n+1}(t)x_{n+1}\|)}2\;dP(t)\\
&\ge \int_{\Omega} \varphi(\|S_n(\omega)\|)\; dP(t) =
\varphi(\|S_n(\omega)\|),
\end{align*} which completes the proof.
\end{proof}

\begin{prop}\label{propmonotone}
Let $\varphi$ be a convex increasing function on $[0,\infty)$ and
$x, y\in X$. Then the function $\psi$ on $\mathbb{R}$ defined by
\[ \psi(\lambda) = \mathbb{E}[\varphi(\|x + \lambda r_1 y\|)]\]  is an
increasing convex function on $[0,\infty)$ with $\psi(\lambda) =
\psi(|\lambda|)$ for every $\lambda\in \mathbb{R}$.
\end{prop}
\begin{proof}By the convexity of $\varphi$, we get
$\psi(\lambda) = \frac12[\psi(\lambda) + \psi(-\lambda) ] \ge
\psi(0)$ for every real $\lambda$, which implies that
$\lambda\mapsto \psi(\lambda)$ is increasing on $[0,\infty)$.
Clearly $\lambda\mapsto \psi(\lambda)$ is  an even function on
$\mathbb{R}$ since $r_1$ is symmetric.
 \end{proof}
Since random variables $\{r_i\}_i$ are independent,
Proposition~\ref{propmonotone} shows that for any two real sequences
$\{\lambda_i\}_{i=1}^n$ and $\{\xi_i\}_{i=1}^n$ satisfying
$|\lambda_i| \le |\xi_i|$ for $i=1, \cdots, n$, and for any $x_0,
x_1, \ldots, x_n$ in $X$, we get
\[ \mathbb{E}[ \varphi\|x_0+ \lambda_1 r_1x_1 + \cdots + \lambda_n r_n x_n\|] \le
\mathbb{E}[\varphi\|x_0+\xi_1 r_1 x_1 +  \cdots + \xi_n r_nx_n \|
].\]

A convex function $\varphi:[0,\infty) \to \mathbb{R}$ is said to be
{\it strictly convex} if
\[\varphi(\frac{s+t}2)<\frac{\varphi(s)+\varphi(t)}2\] holds
for all distinct positive numbers $s,t$. Notice that if  $\varphi$
is strictly convex and increasing on $[0,\infty)$ and $a, b$ are
real numbers with $b\neq0$, then
\[ \mathbb{E}[\varphi(|a + r_2 b|)] =
\mathbb{E}[ \frac{\varphi(|a + r_2 b|) + \varphi(|a - r_2 b|)}2]  >
\varphi(|a|). \]

Now we state our main theorem concerning the randomized series.
\begin{thm}\label{thmbaisics}
Let $\{x_i\}_{i=1}^n$ be a finite sequence in a Banach space $X$,
$\varphi$ a strictly convex, increasing function with $\varphi(0)=0$
and $\{r_i\}_{i=1}^\infty$ symmetric independent random variables
with $\|r_i\|_\infty =1$, $(i=1, 2, \cdots)$.

Suppose that there is a constant $\rho>0$ such that the following
holds: \[ \sup_{\epsilon_1 =\pm 1, \ldots, \epsilon_n =\pm 1} \|
\epsilon_1 x_1 + \cdots + \epsilon_n x_n \|\ge \|x_1\| + \rho.\]
Then there is a constant $\delta=\delta(\rho) >0$ such that
\[ \mathbb{E}[ \varphi(\|x_1 + r_2x_2 + \cdots + r_n x_n \|) ] \ge \varphi(\|x_1\|) +
\delta.\] In particular, if we take $\rho_1  = \min\{\rho, 1/2\}$,
then
\[ \delta = \min\left\{ \varphi(\frac{\rho_1}3)\prod_{i=2}^nP\{|r_i-1|< \frac{\rho_1}{3n}\},
\min_{2\le j\le n}\mathbb{E}\left[\varphi\left(\|x_1\|+r_j
\frac{\rho_1}{3n}\right)\right]-\varphi(\|x_1\|)\right\}.\]
\end{thm}
\begin{proof}We adapt the argument in the proof of Proposition~2.2 in \cite{D}.
We assume that there exist $0<\rho<1/2$ and  signs $\epsilon_2,
\ldots, \epsilon_n$ such that
\[ \|x_1 + \epsilon_2 x_2 + \cdots +\epsilon_n x_n \| \ge \|x_1\|+\rho.\]
Select a unit element $x^*$ in $X^*$ such that $x^*(x_1) =\|x_1\|$
and let $\lambda_i = x^*(x_i)$ for $1\le i\le n$. Now we shall
consider two cases according to the size of $|\lambda_i|$. In the
first case we suppose that $\max_{2\le i\le n} |\lambda_i|\le
\frac{\rho}{3n}$. If $|\eta_i - \epsilon_i| \le \frac{\rho}{3n}$,
$(2\le i \le n)$, then
\[ \|x_1 +\eta_2x_2 + \cdots + \eta_n x_n \| \ge
\|x_1 +\epsilon_2 x_2+\cdots \epsilon_n x_n\| -\frac{\rho}3\ge
\|x_1\|+\frac{2\rho}3.\] Since $|\lambda_1+\eta_2\lambda_2 + \cdots
+ \eta_n\lambda_n | \le |\lambda_1| + \frac\rho3 =\|x_1\| +
\frac{\rho}3$, we get
\[\|x_1 +\eta_2x_2 + \cdots + \eta_n x_n \| \ge
|\lambda_1 +\eta_2 \lambda_2+\cdots \eta_n
\lambda_n|+\frac{\rho}3,\] if $|\eta_i - \epsilon_i| \le
\frac{\rho}{3n}$, $(2\le i \le n)$. Let
\[ F = \bigcap_{j=2}^n \left\{w\in \Omega: |r_j(w) - \epsilon_j|< \frac{\rho}{3n}\right\}\]
and take $T_n = x_1 + r_2x_2 + \cdots + r_n x_n$. Then we have \[
\mathbb{E}[\varphi(\|T_n\|)] = \mathbb{E}[\varphi(\|T_n\|)\chi_F]  +
\mathbb{E}[\varphi(\|T_n\|)\chi_{F^c}].\]  Since $\varphi(a+b) \ge
\varphi(a) + \varphi(b)$, $(a,b\ge 0)$,
\begin{align*}\mathbb{E}[\varphi(\|T_n\|)\chi_F] &\ge
\mathbb{E}[\varphi(|\lambda_1 + r_2\lambda_2 + \cdots + r_n
\lambda_n |+\rho/3)\chi_F]\\&\ge \mathbb{E}[\varphi(|\lambda_1 +
r_2\lambda_2 + \cdots + r_n \lambda_n |)\chi_F] +
\varphi(\rho/3)P(F).\end{align*} Hence
\begin{align*}\mathbb{E}[\varphi(\|T_n\|)]& \ge
\mathbb{E}[\varphi(|\lambda_1 + r_2\lambda_2 + \cdots + r_n
\lambda_n |)] + \varphi(\rho/3)P(F)\\&\ge \varphi(\lambda_1)+
\varphi(\rho/3)P(F) =\varphi(\|x_1\|) +
\varphi(\rho/3)P(F).\end{align*} Notice that
$P(F)=\prod_{i=2}^nP\{w\in \Omega: |r_i(w) - 1|<
\frac{\rho}{3n}\}>0$ for $r_i$'s are independent symmetric random
variables with $\|r_i\|_\infty=1$, $(i=1,2,\cdots)$. In the second
case we suppose that there exists $i_0$ $(2\le i_0\le n)$ such that
$|\lambda_{i_0}|\ge \frac{\rho}{3n}$. It follows from
Proposition~\ref{propsubmartingale}, ~\ref{propmonotone} and strict
convexity of $\varphi$ that
\begin{align*} \mathbb{E}[\varphi(\|T_n \|)] &\ge
\mathbb{E}[\varphi(|\lambda_1 + r_2 \lambda_2 +\cdots + r_n
\lambda_n |)]
\\&\ge \mathbb{E}[\varphi(|\lambda_1+ r_{i_0}\lambda_{i_0}|)]
\\&\ge \mathbb{E}\left[\varphi\left(\left|\lambda_1+
r_{i_0}\frac{\rho}{3n}\right|\right)\right]\\&>\varphi(\lambda_1)=
\varphi(\|x_1\|).\end{align*} The proof is complete.\end{proof}

\section{Representability of $\ell^n_\infty$}

A Banach space $Y$ is said to be {\it representable} in $X$ if, for
each $\lambda>1$, there is a bounded linear map $T: Y \to X$ such
that $\|x\| \le \|Tx\|\le \lambda \|x\|$ for every $x\in Y$. A
Banach space $Y$ is said to be {\it finitely representable} in $X$
if every finite dimensional subspace of $Y$ is representable in $X$.

It is well-known due to the work of B. Mauray and G. Pisier
\cite{MP} that $c_0$ is finitely representable in $X$ if and only if
$c_0$ is finitely representable in $L_p(X)$ for all $1\le p
<\infty$. S. J. Dilworth considered the quatitative version of this
theorem  in \cite{D}, where he showed that if $X$ is a complex
Banach space, $n\ge 2$ and $0<p<\infty$, then
$\ell_\infty^n(\mathbb{C})$ is representable in $X$ if and only if
it is representable in $L_p(X)$. As we see in the next example it is
not true for real Banach spaces.

\begin{example}
Let $X$ be a nontrivial real Banach space. Then $\ell_\infty^2$ is
representable in $L_1([0,1];X)$. Indeed, if we choose the Rademacher
sequence $\{r_n = {\rm sign}( \sin2^n\pi t)\}_{n=1}^\infty$ in
$L_1[0,1]$ and $x_0 \in S_X$, then $x= r_1 x_0$ and $y = r_2 x_0$
are the elements of unit sphere of $L_1(X)$ and they satisfy
\[ \|x+y\|_{L_1(X)} = \|x- y\|_{L_1(X)}=1,\] which means that
$\ell_\infty^2$ is representable in $L_1(X)$.
\end{example}

 The subharmonicity of absolute value of holomorphic functions on
$\mathbb{C}$ plays the crucial role in the proof in \cite{D}. In
this paper, the strict convexity of $\varphi(t) = |t|^p$  $( 1< p
<\infty)$ on $\mathbb{R}$ plays the analogous role, and thus it is
shown here that $\ell_\infty^n$ is representable in $X$ if and only
if it is representable in $L_p(X)$ for every $1<p<\infty$.

The following proposition is a real version of Proposition~2.2 in
\cite{D}.

\begin{prop}\label{criterion-rep}
Let $\{r_i\}_{i=1}^\infty$  be symmetric independent random
variables with $\|r_i\|_\infty =1$, $(i=1, 2, \cdots)$. Suppose that
$X$ is a real Banach space and that $n\ge 2$. The following
properties are equivalent:
\begin{enumerate}
\item $\ell_\infty^n$ is  not  representable in $X$.

\item There exists $\rho>0$ such that whenever $x_1,\ldots, x_n$ are
unit vectors in $X$, then there exist signs $\epsilon_1, \ldots,
\epsilon_n$ such that
\[ \|\epsilon_1 x_1 + \cdots + \epsilon_n x_n \| \ge 1+ \rho.\]

\item There exist strictly convex, increasing function  $\varphi$ on $[0,\infty)$
with $\varphi(0)=0$ and $\rho>0$ such that whenever $x_1, \ldots, x_n$ are unit vectors
in $X$ then \[ \mathbb{E}[ \varphi(\|x_1 + r_2x_2\cdots + r_n
x_n\|)] \ge \varphi(1)+ \rho.\]

\item For each strictly convex,  increasing function $\varphi$ on $[0,\infty)$ with $\varphi(0)=0$, there
is $\rho>0$ such that whenever $x_1, \ldots, x_n$ are unit vectors
in $X$ then \[ \mathbb{E}[ \varphi(\|x_1 + r_2x_2\cdots + r_n
x_n\|)] \ge \varphi(1)+ \rho.\]

\end{enumerate}
\end{prop}
\begin{proof}
The implications $(4)\Rightarrow(3)\Rightarrow(2)\Rightarrow(1)$ are
clear. To show $(1)\Rightarrow (2)$, suppose that $(2)$ fails, so
for any $\rho>0$ there exist unit vectors $x_1, \ldots, x_n$ in $X$
such that for all signs $\epsilon_1, \ldots, \epsilon_n$, we have
\[ \|\epsilon_1 x_1 + \cdots + \epsilon_n x_n \| < 1+ \rho.\] It
follows that for all $\lambda_1, \ldots, \lambda_n$ with
$1=|\lambda_{i_0}| = \max_{1\le i\le n} |\lambda_i|$, we have
\[ \|\lambda_1 x_1 + \cdots + \lambda_n x_n \| < 1+\rho.\]
Hence
\begin{align*} 1= \|\lambda_{i_0}x_{i_0}\| &\le \frac12 \|\lambda_1 x_1 +
\cdots + \lambda_n x_n \| + \frac12 \|\lambda_1x_1 + \cdots +
\lambda_n x_n -2\lambda_{i_0}x_{i_0}\|\\ &\le \frac12 \|\lambda_1x_1
+ \cdots + \lambda_n x_n \| + \frac12(1+\rho).
\end{align*}
Thus, $\|\lambda_1x_1 +\cdots + \lambda_n x_n\|\ge 1-\rho$. Since
$\rho$ is arbitrary, it follows that $\ell_\infty^n$ is
representable in $X$. Now we have only to show that $(2)\Rightarrow
(4)$. Suppose that $(2)$ holds and that $\varphi$ is a strictly
convex, increasing function on $[0, \infty)$ with $\varphi(0)=0$.
There is $0<\rho<1/2$ such that whenever $x_1, \ldots, x_n$ are unit
vectors in $X$, there exist signs $\epsilon_2, \ldots, \epsilon_n$
such that
\[ \|x_1 + \epsilon_2 x_2 + \cdots \epsilon_n x_n \| \ge 1+\rho.\]
Then Theorem~\ref{thmbaisics} shows that $(2)$ implies $(4)$.
\end{proof}

Notice that in the case of the Rademacher sequence $\{r_n\}$, for
every $x_1, \ldots, x_n$ in $X$,
\[\mathbb{E}[\varphi(\|r_1 x_1 + \cdots + r_n x_n \|)] =
\mathbb{E}[\varphi(\|x_1 + r_2 x_2 + \cdots + r_n x_n\|)].\] The
following theorem shows the lifting property of representability of
$\ell_\infty^n$.

\begin{thm}
Suppose $X$ is a Banach space, $(M, \mathfrak{M}, \mu)$ is a measure
space with a measurable subset $A$ satisfying $0<\mu(A)<\infty$, and
$1<p<\infty$, $n\ge2$. Then $\ell_\infty^n$ is representable in $X$
if and only if it is representable in $L_p(M, \mathfrak{M}, \mu;
X)$.
\end{thm}
\begin{proof}
One implication is clear.  To prove the other implication, suppose
that $\ell_\infty^n$ is not representable in $X$ and let $\{r_n\}$
be the Rademacher sequence. By Proposition~\ref{criterion-rep},
there exits $0<\rho<1/2$ such that
 whenever $x_1, \ldots, x_n$ are unit vectors in $X$, we have
\[ \mathbb{E}\|r_1x+r_2 x_2 + \cdots + r_n x_n \|^p \ge (1+\rho)^p.\]
Suppose that $f_1, \ldots, f_n$ are unit vectors in $L_p(X)$. We
define the following functions on $M$. For $w\in M$, let
\begin{align*} q(w)^p &=
\mathbb{E}\|r_1f_1(w) +   \cdots + r_n f_n(w)\|^p,\\
M(w) &= \max\{ \|f_i(w) \| : 1\le i\le n\},\\
m(w) &= \min\{ \|f_i(w) \| : 1\le i\le n\}.
\end{align*}
By Proposition~\ref{propsubmartingale}, $q(w) \ge M(w)$ for all
$w\in M$. Now the argument divides into two cases according to the
relative sizes of $M(w)$ and $m(w)$. In the first case we suppose
that $(1-\rho/3)M(w)\ge m(w)$. Then
\begin{align*}
\frac1n \sum_{i=1}^n \|f_i(w)\|^p &\le \frac{n-1}n M(w)^p + \frac 1n
m(w)^p\\& \le \left( 1- \frac \rho{3n} \right) M(w)^p\\&\le q(w)^p -
\frac \rho{3n} \left( \frac 1n \sum_{i=1}^n
\|f_i(w)\|^p\right)\end{align*} and so
\[ q(w)^p \ge \left( 1+ \frac \rho{3n}\right) \frac{1}n \sum_{i=1}^n
\|f_i(w)\|^p.\] In the second case, we suppose that $(1-\rho/3)M(w)<
m(w)$. Then
\begin{align*}
q(w)^p &\ge (1+\rho)^p m^p(w)\\&\ge
(1+\rho)^p\left(1-\frac\rho3\right)^p \frac 1n \sum_{i=1}^n
\|f_i(w)\|^p\\&\ge \left(1+\frac\rho2\right)^p\frac 1n \sum_{i=1}^n
\|f_i(w)\|^p.
\end{align*} Hence by the Fubini theorem,
\[\mathbb{E}\|rf_1+\cdots+r_nf_n\|_{L^P(X)}^p=\int_{M} q(w)^p\
d\mu \ge \min\left\{\left(1+\frac \rho2\right)^p, \left(1+
\frac\rho{3n}\right)\right\},\] which shows that $\ell_\infty^n$ is
not representable in $L_p(X)$ by Proposition~\ref{criterion-rep}.
The proof is completed.
\end{proof}

\section{Applications to the convexity of Banach spaces}

Recall that  a point $x$ in $S_X$ is an {\it extreme point} of $B_X$
if $\max\{\|x+y\|, \|x-y\|\} =1$ for some $y\in X$ implies $y=0$. A
point $x\in S_X$ is called a {\it strongly extreme point} of $B_X$
if, given $\epsilon>0$, there is a $\delta= \delta(x,\epsilon)>0$
such that
\[ \inf\{ \max\{\|x+y\|, \|x-y\|\} : \|y\|\ge \epsilon\} \ge 1 + \delta.\] A Banach space is said
to be {\it strictly convex} (resp. {\it midpoint locally uniformly
convex}) if every point of $S_X$ is (resp. strongly) extreme point
of $B_X$. A Banach space is called {\it uniformly convex} if, given
$\epsilon>0$, there is a $\delta=\delta(\epsilon)>0$ such that
\[ \inf\{ \max\{\|x+y\|, \|x-y\|\} : \|y\|\ge \epsilon, \|x\|=1\} \ge 1 +
\delta.\]

Theorem~\ref{thmbaisics} gives the following criteria for the
various convexity properties.

\begin{thm}\label{convexity}
Let $X$ be a real Banach space and $\varphi$, a strictly convex
increasing function on $[0, \infty)$ with $\varphi(0)=0$ and $r$, a
symmetric random variable with $\|r\|_\infty=1$. Then
\begin{enumerate}
\item A point $x$ in $S_X$ is an extreme point of $B_X$ if and only if
$\mathbb{E}[\varphi(\|x+ r y\|)] =\varphi(1)$ for $y\in X$ implies
$y=0$.

\item A point $x$ in $S_X$ is  a strongly extreme point of $B_X$ if
and only if for every $\epsilon>0$ there is a $\delta>0$ such that
whenever $\|y\|\ge \epsilon$, we get
\[\mathbb{E}[\varphi(\|x+
r y\|)] \ge \varphi(1) +\delta.\]

\item\label{3uniformconvex} $X$ is uniformly convex if and only if the modulus
$\delta_\varphi(\epsilon)>0$ for every $\epsilon>0$, where $
\delta_\varphi(\epsilon) = \inf\{ \mathbb{E}[\varphi(\|x+ r
y\|)]-\varphi(1) : x\in S_X, \|y\|\ge \epsilon \}.$
\end{enumerate}
\end{thm}
\begin{proof}
We prove only (\ref{3uniformconvex}) because the proof of the others
are similar. Suppose that $X$ is uniformly convex. Given
$\epsilon>0$, there is a $\rho>0$ such that for any $x\in S_X$ and
$y\in X$ with $\|y\|\ge \epsilon$, we have
\[ \max\{\|x+y\|, \|x-y\|\}\ge 1 + \rho.\] Then
Theorem~\ref{thmbaisics} shows that there is $\delta>0$ such that
for any $x\in S_X$ and $y$ with $\|y\|\ge \epsilon$,
\[ \mathbb{E}[\varphi(\|x+r y\|)]\ge\varphi(1) + \delta.\]
Conversely, suppose that $\delta_\varphi(\epsilon)>0$ for every
$\epsilon>0$. Then given $\epsilon>0$ for any $x\in S_X$ and $y\in
X$ with $\|y\|\ge \epsilon$, we have
\[ \max\{\varphi(\|x+y\|), \varphi(\|x-y\|)\} \ge \max_{-1\le t\le 1}
\{ \varphi(\|x+ty\|)\} \ge \mathbb{E}[\varphi(\|x+r y\|)].\] Since
$\varphi$ is strictly increasing we get, for any $x\in S_X$ and
$y\in X$ with $\|y\|\ge \epsilon$,
\[ \max\{\|x+y\|, \|x-y\|\}\ge
\varphi^{-1}(1+\delta_\varphi(\epsilon))>1.\] Therefore $X$ is
uniformly convex and this completes the proof.
\end{proof}

It is worthwhile to notice that Theorem~\ref{convexity} does not
hold if we consider the general increasing convex function $\varphi$
with $\varphi(0)=0$. Indeed, it is easily checked that if
$\varphi(t)=|t|$ and $\{r_n\}_n$ is the Rademacher sequence then for
any nontrivial Banach space $X$,
\[ \mathbb{E}\|x + r_1 x\|=1 \ \ \ (x\in S_X).\]
Consequently, we cannot characterize the extreme point of $B_X$ with
$\varphi(t)=t$.

We shall discuss the uniform convexity of $p$-convexification
$E^{(p)}$ for uniformly monotone Banach lattice $E$. For more
details on Banach lattices, order continuity
 and K\"othe function spaces, see \cite{LT}. For the definition of
$p$-convexification $E^{(p)}$ of $E$ and the addition $\oplus$ and
multiplication $\odot$ there, see \cite{Lee2,LT}. A Banach lattice
is said to be {\it uniformly monotone} (resp. {\it upper locally
uniformly monotone}) if given $\epsilon>0$
\[ M_p(\epsilon) =\inf \{ \|(|x|^p + |y|^p)^{1/p}\|-1: \|y\|\ge
\epsilon, \|x\|=1\}>0\]
\[(\text{resp.}\ \ \ \  N_p(\epsilon;x)=\inf\{\| (|x|^p +
|y|^p)^{1/p} \|-1:\|y\|\ge \epsilon\}>0\ \ \ )\] for some $1\le
p<\infty$.
 It is shown in \cite{Lee2,
Lee} that, given $\epsilon>0$ and $1\le p<\infty$ there is a $C_p>0$
such that for every $\epsilon>0$,
\begin{equation}\label{eqlocalmonotone1}C_p^{-1}
M_1(C_p^{-1}\epsilon^p) \le M_p(\epsilon) \le
M_1(\epsilon).\end{equation}

In the case when $E$ is an order continuous Banach lattice or a
K\"othe function space, we also get the following relations by
 Lemma~2.3 in \cite{Lee}: There is a $C_p>0$ such that every $x\in S_X$ and $\epsilon>0$,
\begin{equation}\label{eqlocalmonotone2}C_p^{-1}
N_1(C_p^{-1}\epsilon^p;x) \le N_p(\epsilon;x) \le
N_1(\epsilon;x).\end{equation}

Notice that relations (\ref{eqlocalmonotone1}),
(\ref{eqlocalmonotone2}) show that if a Banach lattice $E$ is
uniformly monotone then $E^{(1/p)}$ is uniformly monotone
quasi-Banach lattice for $1<p<\infty$. Similarly, if $E$ is upper
locally uniformly monotone order continuous Banach lattice or
K\"othe function space, then $E^{(1/p)}$ is also upper locally
uniformly monotone quasi-Banach lattice for $1<p<\infty$ (cf.
\cite{Lee2}). The characterizations of local uniform monotonicity of
various function spaces have been discussed in \cite{FK}.

In \cite{HKM}, H. Hudzik, A. Kami\'nska and M. Masty\l o showed that
if a K\"othe function space $E$ is uniformly monotone then its
$p$-convexification $E^{(p)}$ is uniformly convex for $1<p<\infty$.
A partial generalization of this result has been studied by the
author in \cite{Lee2}, where it was shown that if a Banach lattice
is uniformly monotone then $E^{(p)}$ is uniformly convex for all
$2\le p<\infty$. In the next theorem, the gap is completed.

\begin{thm}\label{thmmonotonetoconvex}
Let $E$ be a Banach lattice. The following statements are
equivalent.
\begin{enumerate}
\item\label{item1} $E$ is uniformly monotone.

\item\label{item2} $E^{(p)}$ is uniformly convex for all $1<p<\infty$.

\item\label{item3} $E^{(p)}$ is uniformly convex for some $1<p<\infty$.
\end{enumerate}
\end{thm}
\begin{proof}
Proposition~4.4 in \cite{Lee2} shows that uniformly convex Banach
lattice is uniformly monotone. So if we assume (\ref{item3}), then
$E^{(p)}$ is uniformly monotone and $E$ is uniformly monotone. Hence
(\ref{item3})~$\Rightarrow$~(\ref{item1}) is proved. The implication
(\ref{item2})~$\Rightarrow$~(\ref{item3}) is clear. So we have only
to show that (\ref{item1}) implies (\ref{item2}).

We shall use Theorem~\ref{convexity}~(3) with the Rademacher
function $|r|=1$. Let $\epsilon>0$ and let $f, g\in E^{(p)}$ with
$\|f\|_{E^{(p)}}=\|f\|^{1/p}_E= 1$ and $\|g\|_{E^{(p)}}=
\|g\|_E^{1/p}\ge \epsilon$. Recall the following well-known
inequality (cf. Lemma~4.1 \cite{Lee2}) : for any $1<p<\infty$ there
is $C=C(p)$ such that for any reals $s,t$,
\[ \left( \left| \frac{s-t}{C}\right|^2 + \left|
\frac{s+t}{2}\right|^2\right)^{\frac 12} \le \left(\frac{|s|^p+
|t|^p}2\right)^{\frac 1p}.\] Then applying the Krivine functional
calculus to the inequality above, we get
\begin{align}\label{eqlatticeconvex1} \mathbb{E}\left[\|f \oplus (r\odot
g)\|_{E^{(p)}}^p \right]&= \mathbb{E}[\|\ |f^{1/p}+ r g^{1/p}|^p
\|_{E}]\\&= \frac{\mathbb{E}\|\ |f^{1/p}+rg^{1/p}|^p\|_E +
\mathbb{E}\|\ |f^{1/p}-rg^{1/p}|^p\|_E}2\nonumber\\
&\ge \mathbb{E}\left[ \left\|\frac{|f^{1/p}+rg^{1/p}|^p +
|f^{1/p}-rg^{1/p}|^p}2\right\|_{E}\right]\nonumber\\& \ge
\mathbb{E}\left\|\left(|f|^{2/p} + \frac
{|2g|^{2/p}}{C^{2/p}}\right)^{p/2}\right\|_{E}\nonumber\\&=
\left\|\left(|f|^{2/p} + \frac
{|2g|^{2/p}}{C^{2/p}}\right)^{p/2}\right\|_{E}\nonumber.
\end{align}

By (\ref{eqlatticeconvex1}), if $1<p\le 2$, then
$\mathbb{E}\left[\|f \oplus (r\odot g)\|_{E^{(p)}}^p \right]\ge 1 +
M_{2/p}(2\epsilon^p/C)$. In the case of $2<p<\infty$,
(\ref{eqlatticeconvex1}) shows that
\begin{align*}\mathbb{E}\left[\|f \oplus (r\odot g)\|_{E^{(p)}}^p
\right] &\ge  \left\|\left(|f|^{2/p} + \frac
{|2g|^{2/p}}{C^{2/p}}\right)^{p/2}\right\|_{E}\\
&\ge \left\|\left(|f| + \frac {|2g|}{C}\right)\right\|_{E}\ge 1+
M_1(2\epsilon^p/C).
\end{align*}
Hence
\[ \mathbb{E}[ \|f\oplus (r\odot g) \|_{E^{(p)}}^p] \ge 1+M_{\max\{1, 2/p  \}}(2\epsilon^p/C)\]
completes the proof.
\end{proof}

Now we discuss the the local version of
Theorem~\ref{thmmonotonetoconvex}. A point $x\in S_X$ in a complex
Banach space $X$ is said to be a {\it complex strongly extreme
point} if  there is $0<p<\infty$ such that given $\epsilon>0$,
\[ H_p(\epsilon; x) =\inf \left\{\left( \frac{1}{2\pi} \int_0^{2\pi} \|x+
e^{i\theta}y\|^p\ d\theta\right)^{1/p}-1: \|y\|\ge \epsilon \right\}
> 0.\]

It is known in \cite{DHM} that $x\in S_X$ is a complex strongly
extreme point if and only if for every $\epsilon>0$,
\[ H_\infty(\epsilon; x) = \inf\{ \max_{0\le \theta\le
2\pi}\|x+e^{i\theta}y\|-1: \|y\|\ge \epsilon \}>0.\] For more
details about these moduli, see \cite{DGT, D, DHM}. A complex Banach
space $X$ is said to be {\it locally uniformly complex convex} if
every point of $S_X$ is a complex strongly extreme point.

\begin{thm}
Let $E$ be an order continuous Banach lattice or a K\"othe function
space. Then the following are equivalent:
\begin{enumerate}
\item $E$ is upper locally uniformly monotone.

\item $E^{\mathbb{C}}$ is locally uniformly complex convex.

\item $E^{(p)}$ is midpoint locally uniformly convex for all
$1<p<\infty$.

\item $E^{(p)}$ is midpoint locally uniformly convex for some
$1<p<\infty$.

\end{enumerate}
\end{thm}
\begin{proof}First we prove the equivalence of
(1) and (2). We shall use a similar argument as in the proof of
\cite[Proposition~3.7]{CHL} in the sequence space. Suppose that $E$
is locally uniformly complex convex. Then for each $x\in S_E$ and
$\epsilon>0$ there is $\delta=\delta(x,\epsilon)>0$ such that for
all $y \in X$ with $\|y\|\ge \epsilon$
\[ \|\ |x| + |y|\  \|\ge \frac 1{2\pi} \int_0^{2\pi} \|x+ e^{i\theta} y\| \; d\theta \ge 1 + \delta.\]
So $X$ is upper locally uniformly monotone.

 For the converse, suppose that $E$ is upper locally uniformly
monotone. Now, if we use Theorem~7.1 in \cite{DGT}, then we have for
every pair $x,y$ in $E$,
\[ \frac{1}{2\pi}\int_0^{2\pi}\|x+e^{i\theta} y\|\; d\theta
\ge \left\| \left(|x|^2 + \frac12|y|^2\right)^{1/2}\right\|. \]
Hence for every $x\in S_X$ and $\epsilon>0$, we get
\[ \frac{1}{2\pi}\int_0^{2\pi}\|x+e^{i\theta} y\|\; d\theta
\ge 1+ N_2(\epsilon/\sqrt{2};x).\] Therefore, the upper local
uniform monotonicity of $E$ implies the local uniform complex
convexity of $E$.

For (1)$\Rightarrow$(3), fix $f$ with $\|f\|_{E^{(p)}}=
\|f\|^{1/p}=1$ and for any $g\in E^{(p)}$ with $\|g\|_{E^{(p)}}=
\|g\|^{1/p}\ge\epsilon$, (\ref{eqlatticeconvex1}) holds. Hence
\[ \mathbb{E}[ \|f\oplus (r\odot g) \|_{E^{(p)}}^p] \ge 1+N_{\max\{1, 2/p  \}}(2\epsilon^p/C;f),\]
which shows that  (1)$\Rightarrow$(3) holds.

The implication (3)$\Rightarrow$(4) is clear. Finally assume that
(4) holds. Note that every midpoint locally uniformly convex Banach
lattice is upper locally uniformly monotone. Indeed, if $x\in S_X$
and $\epsilon$, there is $\delta>0$ such that
\[ 1+\delta \le \max\{ \|x+y\|, \|x-y\|\} \le \| \ |x| + |y|\ \|.\]
Since the midpoint local uniform convexity of $E^{(p)}$ implies the
upper local uniform monotonicity of $E^{(p)}$, $E$ is upper locally
uniformly monotone. This completes the proof.\end{proof}

Let $X$ be a real Banach space and $\Delta$ be the open unit disk in
$\mathbb{C}$. Let $(f_n)$ be a sequence of continuous functions from
$\Delta$ into $X$ and $f:\Delta\to X$ be continuous. We say that
$(f_n)$ converges to $f$ with respect to the topology of norm
uniform convergence on compact subsets of $\Delta$ if $\lim_{n\to
\infty} \sup\{\|f_n(z)- f(z)\|: z\in K \}=0$ for all compact subsets
$K$ of $\Delta$. We will denote by $\beta$ the topology of norm
uniform convergence on compact subsets of $\Delta$.

A Banach space $X$ is said to have {\it Kadec-Klee property} with
respect to topology $\tau$ $(KK(\tau))$ if whenever $(x_n)$ is a
sequence in $X$ and $x\in X$ satisfy $\|x_n\|=\|x\|=1$ for all $n\in
\mathbb{N}$ and $\tau$-$\lim_n x_n = x$, then $\lim_n\|x_n - x\|=
0$.

A function $f:\Delta\to X$ is {\it harmonic} if $f$ is twice
continuously differentiable and if the Laplacian of $f$ is zero. It
is known  \cite{H} that $f:\Delta\to X$ is harmonic if and only if
$x^*f$ is harmonic for all $x^*\in X^*$ if and only if there is a
sequence $\{a_n\}_n\subset X$ so that for all $0\le r<1$, $\theta\in
\mathbb{R}$,
\[ f(re^{i\theta}) = \sum_{n=-\infty}^\infty a_n r^{|n|}e^{in \theta}
,\] where the series is absolutely and locally uniformly convergent.

We now define $h^p(\Delta; X)$ for $1<p<\infty$ by
\[ h^p(\Delta; X) = \{ f:\Delta\to X : f\  \text{is harmonic and}\
\|f\|_p<\infty\},\] where
\[ \|f\|_p = \sup_{0\le r<1} \left( \frac{1}{2\pi}\int_0^{2\pi} \|f(re^{i\theta})\|^p
\;d\theta \right)^{1/p}.\] It is easy to see that on $h^p(\Delta;
X)$, $\|\cdot\|_p$ is $\beta$-lower semicontinuous function.

It is shown by P. N. Dowling and C. J. Lennard \cite{DL} that if
$h^p(\Delta; X)$ has $KK(\beta)$, then $X$ is strictly convex and
has the Radon-Nikod\'ym property.

In fact, the following proposition is a consequence of the results
in \cite{DHS}. We present an easy proof.

\begin{thm}
If $h^p(\Delta; X)$ has $KK(\beta)$, then $X$ is midpoint locally
uniformly convex.
\end{thm}
\begin{proof}
Suppose that $X$ is not locally uniformly convex. Then applying
Theorem~\ref{convexity} with $r(\theta) = \cos\theta$, there exist
an $\epsilon>0$, a sequence $(x_n)$ in $X$ and $x\in S_X$ such that
$\|x_n\|\ge \epsilon$ and
\[\lim_{n\to \infty} \frac{1}{2\pi}\int_0^{2\pi}
\|x+(\cos \theta) x_n\|^p \;d\theta= 1.\] Define $f_n: \Delta\to X$
by
\[ f_n(z) = x + \frac 12 (z^n + \bar{z}^n)x_n\]
and $f:\Delta\to X$ by $f(z)=x$. Then it is easy to see that
$\beta$-$\lim_n f_n(z) = f(z)$. Notice that \begin{align*}
\|f_n\|^p_p & =\sup_{0\le r<1} \frac{1}{2\pi}\int_0^{2\pi} \|x +
r^n\cos (n\theta)x_n \|^p \;d\theta \\ &=
\frac{1}{2\pi}\int_0^{2\pi} \|x + \cos\theta x_n\|^p \;d
\theta.\end{align*} Then $\lim_n \|f_n\|_p = 1=\|f\|_p$. However
\begin{align*} \|f_n - f\|_p^p & =\sup_{0\le r<1} \frac{1}{2\pi}\int_0^{2\pi} \|r^n\cos (n\theta)x_n
\|^p \;d \theta \\ &= \frac{\|x_n\|^p}{2\pi} \int_0^{2\pi}
|\cos\theta|^p \; d \theta\ge \frac{\epsilon^p}{2\pi} \int_0^{2\pi}
|\cos\theta|^p \; d \theta.\end{align*} Hence $h^p(\Delta; X)$ fails
to have $KK(\beta)$. The proof is complete.
\end{proof}
It is worthwhile to remark here that it has been shown in \cite{DHS}
that $h^p(\Delta; X)$ has $KK(\beta)$ if and only if $X$ has the
Radon-Nikod\'ym property and every element of $S_X$ is a denting
point of $B_X$, which is called {\it property $(G)$}. It is easy to
see that a Banach space with property $(G)$ is midpoint locally
uniformly convex.

\subsection*{Acknowledgment}
The author thanks A. Kami\'nska for useful comments.

\bibliographystyle{amsplain}

\end{document}